# Liouville-type theorems on conformal mappings and their application

S. E. Stepanov, I. I. Tsyganok

**Abstract.** In the present paper we prove Liouville-type theorems: non-existence theorems for conformal mappings of complete Riemannian manifolds. In addition, we give applications of these results to theory of conharmonic transformations. A part of these results was announced in our reports on the conferences "Differential Geometry and its Applications" (July 11-15, 2016, Brno, Czech Republic) and "International Conference on Algebra, Analysis and Geometry" (June 26-July 2, 2016, Kazan, Russia).

## Introduction

Conformal mappings of smooth manifolds are important concepts in Riemannian geometry. Liouville's theorem made clear already in the 19th century that in dimensions $n \geq 3$ conformal mappings are more rigid than in dimension 2. In addition, more than forty years ago, Yau has showed in [1] the effectiveness of the Liouville-type theorems in the theory of conformal mappings of an $n$-dimensional ($n \geq 3$) complete Riemannian manifold ($M, g$) whose scalar curvature $s$ is nonpositive and sectional curvature $sec$ ($X, Y$) has a lower bounded, i.e. $sec$ ($X, Y$) $\geq - C$ for some positive constant $C$. In the present paper we prove Liouville-type theorems: non-existence theorems for conformal mappings of complete Riemannian manifolds (without any assumption on their sectional curvatures). To do this, we will use a generalization of the Bochner technique (see [2]). In conclusion, we give applications of these results to the theory of conharmonic transformations.

A part of these results was announced in our reports on the conferences "Differential Geometry and its Applications" (July 11-15, 2016, Brno, Czech Republic) and "International Conference on Algebra, Analysis and Geometry" (June 26-July 2, 2016, Kazan, Russia).

# 1. Subharmonic and superharmonic functions

Let $(M, g)$ be an $n$-dimensional $(n \geq 2)$ Riemannian manifold. Recall here that $\varphi \in C^2 M$ is *subharmonic* (resp. *superharmonic* or *harmonic*) if $\Delta \varphi \geq 0$ (resp. $\Delta \varphi \leq 0$ or $\Delta \varphi = 0$) for the Laplace-Beltrami operator $\Delta \varphi = \operatorname{div}(\operatorname{grad} \varphi)$ (see [3, p. 338]).

In particular, if $(M, g)$ is compact then every harmonic (subharmonic and superharmonic) functions is constant by the lemma of Hopf (see [3, p. 338-339]).

We prove the following lemma on superharmonic functions which consists of two statements that are analogies of two Yau propositions on subharmonic functions from [4]. Firstly, Yau has argued in [4, p. 660] that on a complete Riemannian manifold $(M, g)$ each subharmonic function $u \in C^2 M$ whose gradient has integrable norm on $(M, g)$ must be harmonic. Secondly, he has shown in [4, p. 663] that on a complete Riemannian manifold each non-negative subharmonic function $u \in C^2 M$ such that $\int_M u^p dVol_g < \infty$ for some $1 < p < \infty$ must be constant. In particular, if the volume of $(M, g)$ is infinite, then $u = 0$.

**Lemma**. *If $(M, g)$ is a connected complete Riemannian manifold (without boundary), then any superharmonic function $\varphi \in C^2 M$ with $\|\operatorname{grad} \varphi\| \in L^1(M, g)$ is harmonic and each non-positive superharmonic function $\varphi \in C^2 M$ such that $\varphi \in L^p(M, g)$ for some $1 < p < \infty$ must be constant. In particular, if the volume of $(M, g)$ is infinite, then $\varphi = 0$.*

**Proof**. On the one hand, if we assume that $u = -\varphi$ for any superharmonic function $\varphi \in C^2 M$ then the conditions $\Delta \varphi \leq 0$ and $\|\operatorname{grad} \varphi\| \in L^1(M, g)$ which must be satisfy for the superharmonic function $\varphi$ can be written in the form $\Delta u \geq 0$ and $\|\operatorname{grad} u\| \in L^1(M, g)$. In this case, using the Yau statement for subharmonic functions we conclude that $\Delta u = 0$ and hence $\varphi = -u$ is a harmonic function. On the other hand, the function $u = -\varphi$ for any superharmonic function $\varphi \in C^2 M$ which satisfies the conditions $\varphi \leq 0$, $\Delta \varphi \leq 0$ and $\int_M |\varphi|^p dVol_g < \infty$ for some $1 < p < \infty$ must be satisfied the following conditions $u \geq 0$, $\Delta u \geq 0$ and $\int_M u^p dVol_g < \infty$ for some $1 < p < \infty$. Therefore, $u$ is a constant function and hence $\varphi = -u$ is a constant function too. It is obvious that if the volume of $(M, g)$ is infinite, then $\varphi = 0$.

## 2. Conformal diffeomorphisms of complete Riemannian manifolds

Let $(M, g)$ and $(\overline{M}, \overline{g})$ be pseudo-Riemannian or Riemannian manifolds such that $\dim M = \dim \overline{M} = n$ for any $n \geq 3$. Then a diffeomorphism $f : (M, g) \to (\overline{M}, \overline{g})$ is called *conformal* if it preserves angles between any pair curves. In this case, $\overline{g} = e^{2\sigma} g$ for some scalar function $\sigma$ and for a pointwise identity on $M$ (see [5, p. 89]). If $M = \overline{M}$ then $f : (M, g) \to (M, \overline{g})$ we call a *conformal transformation*. In addition, if the function $\sigma$ is a constant then $f$ is a *homothetic mapping* (resp. *homothetic transformation*). In particular, if $\sigma = 0$ then $f$ is an *isometric mapping* (resp. *isometric transformation*).

If $\sigma \in C^2 M$ then for each pair of corresponding points $x \in M$ and $\overline{x} = f(x) \in \overline{M}$ we have the equation (see [5, p. 90])

$$e^{2\sigma} \overline{s} = s - 2(n-1)\Delta \sigma - (n-1)(n-2)\|\operatorname{grad} \sigma\|^2 \qquad (2.1)$$

where $s$ and $\overline{s}$ denote the scalar curvatures of $(M, g)$ and $(\overline{M}, \overline{g})$, respectively. In the case when $(M, g)$ and $(\overline{M}, \overline{g})$ are Riemannian manifolds we can formulate the following Liouville-type non-existence theorem.

**Theorem 1.** *Let $(M, g)$ be an n-dimensional ($n \geq 3$) complete Riemannian manifold and $f : (M, g) \to (\overline{M}, \overline{g})$ be a conformal diffeomorphism onto another Riemannian manifold $(\overline{M}, \overline{g})$ such that $\overline{g} = e^{2\sigma} g$ and $s \leq e^{2\sigma} \overline{s}$ for some function $\sigma \in C^2 M$ and the scalar curvatures $s$ and $\overline{s}$ of $(M, g)$ and $(\overline{M}, \overline{g})$, respectively. Then the following propositions are true.*

1. *If $\|\operatorname{grad} \sigma\| \in L^1(M, g)$, then $f$ is a homothetic mapping.*

2. *If $\sigma$ is non-positive function and $\sigma \in L^p(M, g)$ for some $1 < p < \infty$ then $f$ is a homothetic mapping. In particular, if the volume of $(M, g)$ is infinite, then $f$ is an isometric mapping.*

**Proof.** If $f : (M, g) \to (\overline{M}, \overline{g})$ is a conformal diffeomorphism a connected complete Riemannian manifold $(M, g)$ onto another Riemannian manifold $(\overline{M}, \overline{g})$ such that $\overline{g} = e^{2\sigma} g$ for some function $\sigma \in C^2 M$, then from (2.1) we obtain

$$2(n-1)\Delta \sigma = s - e^{2\sigma} \overline{s} - (n-1)(n-2)\|\operatorname{grad} \sigma\|^2. \qquad (2.2)$$

Let $s \leq e^{2\sigma} \bar{s}$ then (2) shows $\Delta \sigma \leq 0$. It means that $\sigma$ is a superharmonic function. By the condition of our theorem, the gradient of $\sigma$ has integrable norm on $(M, g)$ and we obtain from (4.2) that $\Delta \sigma = 0$ (see our Lemma). In this case, $\sigma$ is a harmonic function. Since $n \geq 3$, we see from (2.2) that $\sigma$ is constant. In the other hand, if $\sigma$ is a non-positive function such that $s \leq e^{2\sigma} \bar{s}$ and $\sigma \in L^p(M, g)$ for some $1 < p < \infty$ then using the Lemma we can conclude that $\sigma$ is a constant function. It is obvious that if the volume of $(M, g)$ is infinite, then $\sigma = 0$ (see our Lemma). The proof of the theorem is complete.

In particular, if we assume that $s \leq 0$ and $\bar{s} \geq 0$ in the condition of our theorem, then the inequality $s \leq e^{2\sigma} \bar{s}$ holds. Then, as a result the proofs of the Theorem 1, we can conclude that $s = \bar{s} = 0$. Therefore we have the following statement which complements Theorem 1 of Yau in [1].

**Corollary 1**. *Let $(M, g)$ be an n-dimensional $(n \geq 3)$ complete Riemannian manifold and $f : (M, g) \to (\bar{M}, \bar{g})$ be a conformal diffeomorphism onto another Riemannian manifold $(\bar{M}, \bar{g})$ such that $\bar{g} = e^{2\sigma} g$ for some function $\sigma \in C^2 M$, $s \leq 0$ and $\bar{s} \geq 0$ for the scalar curvatures $s$ and $\bar{s}$ of $(M, g)$ and $(\bar{M}, \bar{g})$, respectively. If the one of the following conditions holds*:

1. $\|\operatorname{grad} \sigma\| \in L^1(M, g)$,

2. $\sigma \in L^p(M, g)$ for some $1 < p < \infty$ and $\sigma \leq 0$,

*then f is a homothetic mapping and $s = \bar{s} = 0$. If in the second case the volume of $(M, g)$ is infinite, then f is an isometric mapping.*

Let $\sigma = \log \lambda$ for some positive scalar function $\lambda \in C^2 M$ then

$$\Delta \sigma = \lambda^{-1} \Delta \lambda - \lambda^{-2} \|\operatorname{grad} \lambda\|^2, \qquad \|\operatorname{grad} \sigma\|^2 = \lambda^{-2} \|\operatorname{grad} \lambda\|^2.$$

In this case, (2.2) can be rewritten in the following equivalent form

$$2(n-1)\lambda \Delta \lambda = \lambda^2 (s - \lambda^2 \bar{s}) - (n-1)(n-4)\|\operatorname{grad} \lambda\|^2. \tag{2.3}$$

If $s \geq \lambda^2 \bar{s}$ then for the case $n \leq 4$ from (2.3) we obtain that $\lambda \Delta \lambda \geq 0$. On the other hand, Yau has proved in [4, p. 664] that if a smooth function $\lambda \in C^2 M$ on a complete Riemannian manifold (*M, g*) such that $\lambda \Delta \lambda \geq 0$, then either $\int_M |\lambda|^p dV_g = \infty$ for all $p \neq 1$ or $\lambda = $ constant. Therefore, in the case when (*M, g*) and $(\overline{M}, \overline{g})$ are Riemannian manifolds we have

**Theorem 2.** *Let* (*M, g*) *be a n-dimensional* (*n* = 3, 4) *complete Riemannian manifold and* $f:(M,g) \to (\overline{M}, \overline{g})$ *be a conformal diffeomorphism onto another Riemannian manifold* $(\overline{M}, \overline{g})$ *such that* $\overline{g} = \lambda^2 g$ *and* $s \geq \lambda^2 \overline{s}$ *for some positive function* $\lambda \in C^2 M$ *and for the scalar curvatures s and* $\overline{s}$ *of* (*M, g*) *and* $(\overline{M}, \overline{g})$, *respectively. If* $\lambda \in L^p(M,g)$ *for some* $p \neq 1$, *then f is a homothetic mapping.*

In particular, if we assume that $s \geq 0$ and $\overline{s} \leq 0$ in the condition of Theorem 2, then one can verify that in this case *f* is a homothetic mapping and $s = \overline{s} = 0$. Therefore we have the following statement which generalizes Theorem 2 from [1] for the cases *n* = 3, 4.

**Corollary 2.** *Let* (*M, g*) *be a n-dimensional* (*n* = 3, 4) *complete Riemannian manifold and* $f:(M,g) \to (\overline{M}, \overline{g})$ *be a conformal diffeomorphism onto another Riemannian manifold* $(\overline{M}, \overline{g})$ *such that* $\overline{g} = \lambda^2 g$ *for some positive function* $\lambda \in C^2 M$ *and* $\lambda \in L^p(M,g)$ *for some* $p \neq 1$. *If* $s \geq 0$ *and* $\overline{s} \leq 0$ *for the scalar curvatures s and* $\overline{s}$ *of* (*M, g*) *and* $(\overline{M}, \overline{g})$, *respectively, then f is a homothetic mapping and* $s = \overline{s} = 0$.

On the other hand, the well known that (see [3, p. 338])

$$\Delta \lambda^2 = 2\lambda \Delta \lambda + \|grad\, \lambda\|^2.$$

Then using this equation we can rewrite (2.3) in the form

$$(n-1)\Delta \lambda^2 = \lambda^2 (s - \lambda^2 \overline{s}) - (n-1)(n-5)\|grad\, \lambda\|^2. \qquad (2.4)$$

In this case, if we assume that *n* = 5 and $s \geq \lambda^2 \overline{s}$ then from (2.6) we obtain $\Delta \lambda^2 \geq 0$. Therefore, we can formulate the following two statements.

**Theorem 3**. *Let* (*M, g*) *be an five-dimensional complete Riemannian manifold and* $f:(M,g) \to (\overline{M}, \overline{g})$ *be a conformal diffeomorphism onto another Riemannian manifold* $(\overline{M}, \overline{g})$ *such that* $\overline{g} = \lambda^2 g$ *for some positive function* $\lambda \in C^2 M$. *If* $s \geq \lambda^2 \overline{s}$ *and* $\lambda \in L^{2p}(M,g)$ *for some* $1 < p < \infty$, *then f is a homothetic mapping and* $s = \overline{s} = 0$.

**Corollary 3.** *Let (M, g) be a five-dimensional complete Riemannian manifold and $f:(M,g)\to(\overline{M},\overline{g})$ be a conformal diffeomorphism onto another Riemannian manifold $(\overline{M},\overline{g})$ such that $\overline{g}=\lambda^2 g$ for some positive function $\lambda\in C^2 M$. If $s\geq 0$, $\overline{s}\leq 0$ and $\lambda\in L^{2p}(M,g)$ for some $1<p<\infty$, then f is a homothetic mapping and $s=\overline{s}=0$.*

If we assume that $\lambda=u^{\frac{2}{n-2}}$ then (2.3) immediately gives

$$\frac{4(n-1)}{n-2}\Delta u = su - \overline{s}\,u^{\frac{n+2}{n-2}}. \tag{2.5}$$

In the case of the Riemannian manifolds (M, g) and $(\overline{M},\overline{g})$ the equation (2.5) is the classical *Yamabe equation* (see [6, p. 39]). The equation (2.5) we rewrite in the form

$$\frac{4(n-1)}{n-2}\Delta u = u(s-\lambda^2\,\overline{s}). \tag{2.6}$$

Then for $s\geq\lambda^2\overline{s}$ from (2.6) we obtain that $\Delta u\geq 0$. On the other hand, Yau has shows in [4, p. 663] that if u is a non-negative subharmonic function defined on a complete Riemannian manifold (M, g), then $\int_M u^p dV_g = \infty$ for all $p>1$, unless $u=constant$. Therefore, in the case when (M, g) and $(\overline{M},\overline{g})$ are Riemannian manifolds we have the following Liouville-type non-existence theorem.

**Theorem 4.** *Let (M, g) be a n-dimensional (n ≥ 3) complete Riemannian manifold and $f:(M,g)\to(\overline{M},\overline{g})$ be a conformal diffeomorphism onto another Riemannian manifold $(\overline{M},\overline{g})$ such that $\overline{g}=\lambda^2 g$ and $\lambda^{(n-2)/2}\in L^p(M,g)$ for some positive function $\lambda\in C^2 M$ and for some $p\neq 1$. If $s\geq\lambda^2\overline{s}$ for the scalar curvatures s and $\overline{s}$ of (M, g) and $(\overline{M},\overline{g})$, respectively, then f is a homothetic mapping.*

In particular, if we assume that $s\geq 0$ and $\overline{s}\leq 0$ in the condition of Theorem 3, then we can prove that f is a homothetic mapping and $s=\overline{s}=0$. Therefore we have

**Corollary 4.** *Let (M, g) be a n-dimensional (n ≥ 3) complete Riemannian manifold and $f:(M,g)\to(\overline{M},\overline{g})$ be a conformal diffeomorphism onto another Riemannian manifold $(\overline{M},\overline{g})$ such that $\overline{g}=\lambda^2 g$ and $\lambda^{(n-2)/2}\in L^p(M,g)$ for some positive function $\lambda\in C^2 M$*

*for some* $p \neq 1$. *If* $s \geq 0$ *and* $\bar{s} \leq 0$ *for the scalar curvatures s and $\bar{s}$ of* $(M, g)$ *and* $(\overline{M}, \overline{g})$, *respectively, then f is a homothetic mapping and* $s = \bar{s} = 0$.

### 3. An application to the theory of conharmonic transformations

Ishi in [7] called $f:(M,g) \to (M,\bar{g})$ a conharmonic transformation if it is a conformal transformation $\bar{g} = e^{2\sigma}g$ for some scalar function $\sigma \in C^2 M$ satisfying the equation

$$\Delta \sigma = -\frac{n-2}{2} \|\operatorname{grad} \sigma\|^2 \qquad (3.1)$$

for any $n \geq 3$. The conharmonic transformations introduced by Ishi are a subgroup of the group of conformal transformations which preserve the harmonicity of certain class of smooth functions (see [7]).

From (3.1) we conclude that $\sigma$ is a superharmonic function. Then the following corollary is obvious from the Theorem 1.

**Corollary 5**. *Let* $f:(M,g) \to (M,\bar{g})$ *be a conharmonic transformation of an n-dimensional* $(n \geq 3)$ *complete Riemannian manifold* $(M, g)$, *i.e.* $\bar{g} = e^{2\sigma}g$ *for some function* $\sigma \in C^2 M$ *which satisfies the equation (3.1). If $\sigma$ has a gradient with integrable norm on* $(M, g)$, *then the function $\sigma$ is constant and f is a homothetic transformation.*

Let $\sigma = \log \lambda$ for some positive scalar function $\lambda \in C^2 M$ then (3.1) can be rewritten in the following equivalent form

$$2\lambda \Delta \lambda = (n-4)\|\operatorname{grad} \lambda\|^2. \qquad (3.2)$$

In this case, we can formulate a proposition that is an analog of Theorem 3.

**Corollary 6**. *Let* $f:(M,g) \to (M,\bar{g})$ *be a conharmonic transformation of an n-dimensional* $(n \geq 4)$ *complete Riemannian manifold* $(M, g)$, *i.e.* $\bar{g} = \lambda^2 g$ *for some positive function* $\lambda \in C^2 M$ *which satisfies the equation (3.2). If* $\lambda \in L^p(M,g)$ *for some* $p \neq 1$, *then f is a homothetic transformation.*

In particular, if $n = 4$ then from (3.2) we obtain $\Delta \lambda = 0$. On the other hand, Yau proved that there is no non-constant positive harmonic function on a complete Riemannian manifold with non-negative Ricci curvature (see [8; p. 217]). Therefore, we have

**Theorem 5**. Any *conharmonic transformation of a four-dimensional complete Riemannian manifold* (*M, g*) *with non-negative Ricci curvature is homothetic transformation.*

**Remark**. The corollaries 4 and 5 generalize Proposition 4.7 from [9] on conharmonic transformations of compact manifolds.

**Acknowledgements.** Authors are supported by RBRF grant 16-01-00053-a (Russia).

*Author's address:*

Sergey E. Stepanov[1,2]
[1] All Russian Institute for Scientific and Technical Information
of the Russian Academy of Sciences, 20, Usievicha street,
125190 Moscow, Russian Federation
[2] Department of Mathematics,
Finance University under the Government of Russian Federation,
125468 Moscow, Leningradsky Prospect, 49-55, Russian Federation,
E-mail: s.e.stepanov@mail.ru

Irina I. Tsyganok
Department of Mathematics,
Finance University under the Government of Russian Federation,
125468 Moscow, Leningradsky Prospect, 49-55, Russian Federation,
E-mail: i.i.tsyganok@mail.ru